\documentclass[11pt, a4paper]{article}
\usepackage{}
\usepackage{mathrsfs}
\usepackage{}
\usepackage{amsfonts}
\usepackage{amsmath,amssymb,bm,graphicx,latexsym,graphics,epsfig,subfigure,color}
\usepackage{amsthm}
\usepackage{bm}   
\usepackage{amsmath,amsfonts,amssymb} 
\usepackage{bbm}
\usepackage{color}
\usepackage{xcolor}
\usepackage{graphicx}
\usepackage{epstopdf}
\usepackage{epsfig}
\usepackage{setspace}
\usepackage{indentfirst}
\usepackage{booktabs}
\usepackage{float}
\usepackage[utf8]{inputenc}
\usepackage[title]{appendix}
\usepackage{amsthm}
\usepackage[thicklines]{cancel}
\usepackage{mathrsfs}
\usepackage{multirow}
\usepackage{diagbox}
\usepackage[colorlinks=true,linkcolor=blue,anchorcolor=red,urlcolor=blue,citecolor=blue]{hyperref} 
\usepackage{amssymb}  
\usepackage{geometry}
\usepackage[authoryear]{natbib}

\usepackage{hyperref}
\hypersetup{bookmarksnumbered, bookmarksopen,
	colorlinks,citecolor=blue,linkcolor=blue,CJKbookmarks}

\usepackage{hyperref}
\usepackage{url}
\usepackage{algorithm}
\usepackage{algpseudocode}
\usepackage{multirow}
\usepackage{booktabs}  
\usepackage{graphicx}
\usepackage{comment}

\newtheorem{theorem}{Theorem}

\newtheorem{proposition}{Proposition}
\newtheorem{remark}{Remark}
\newtheorem{corollary}{Corollary}
\newtheorem{example}{Example}

\usepackage{bbm}
\usepackage{color}
\usepackage{colortbl}
\definecolor{text1}{cmyk}{1,.65,0,0} 
\definecolor{text2}{rgb}{1,0,0} 
\definecolor{text3}{cmyk}{0,0,0,1} 
\definecolor{text4}{cmyk}{0,0,0,0.5} 
\definecolor{text5}{cmyk}{1.0,0.0,1.0,0} 

\title{Unbiased estimation in one-parameter exponential families for the inverse of the natural parameter with extensions}
\author{Pankaj Bhagwat$^{a}$ \& \'Eric Marchand$^{b}$} 

\begin{document}

\maketitle

\begin{center}
$a$ University of Alberta (e-mail: pbhagwat@ualberta.ca), Department of Mathematical \& Statistical Sciences, Edmonton, AB, Canada T6G 2N8 \\
$b$ Université de Sherbrooke  (e-mail: Eric.Marchand@usherbrooke.ca), Département de mathématiques, Sherbrooke, QC, Canada J1K 2R1
\end{center}

\vspace{5mm}

\begin{center}
    \begin{abstract}
For one-parameter continuous exponential families, we identify an unbiased estimator of the inverse of the natural parameter $\theta$ for cases where $\theta > 0$, extending an earlier result of \cite{voinov1985unbiased} applicable to a normal model.   We provide various applications for Gamma models, Inverse Gaussian models, distributions obtained by truncation, and ratios of normal means.  Moreover, we extend the findings to estimating negative powers $\theta^{-k}$, and more generally to complete monotone functions $q(\theta)$.
\end{abstract}
\end{center}

\vspace{5mm}



\section{Introduction}

\noindent   Unbiased estimation is a unifying theme in several facets of statistical theory and practice, for parametric and non-parametric models, and also for non-asymptotic and asymptotic settings. Methods and theoretical foundations for finding optimal unbiased estimators, when they exist, are well established, in linear models or under the criteria of uniformly minimum variance unbiased (UMVU) estimation, for instance (e.g., \cite{lehmann1998point}, and the references therein).   The study of efficiency for unbiased or nearly unbiased estimators is also highly relevant (e.g., \cite{doss1989price}). In this paper, we address the issue of the very existence of an unbiased estimator and how this can be impacted by the parameter space setting. \\

\noindent For one-parameter exponential families, we record situations where an unbiased estimator exists for the inverse $\gamma(\theta) = 1/\theta$ of a natural parameter $\theta$; which is sometimes achievable by virtue of a lower bound parametric restriction on $\theta$, in contrast to the unrestricted parameter setting.   As an exemplar, consider a normal model $X \sim N(\theta, \sigma^2)$ with known $\sigma^2$.   For the natural parameter space $\theta \in \mathbb{R}$, it is known and verifiable that there does not exist an unbiased estimator of $1/\theta$ (see Theorem \ref{theoremnon-existence2}).   However, as demonstrated by \cite{voinov1985unbiased}, there {\bf does exist} an unbiased estimator $\delta_0(X)$ when $\theta>0$, i.e.,
\begin{equation}
\nonumber
\mathbb{E}_{\theta}\big(\delta_0(X)  \big) \, = \, \frac{1}{\theta}\,, \hbox{ for all } \theta > 0\,.
\end{equation}

\noindent We revisit this result and, furthermore, we show that a similar result holds for one-parameter continuous exponential families, with extensions to parametric functions of the form $\gamma(\theta)\,=\,1/\theta^k$ with $k>0$ and more generally to completely monotone parametric functions $q(\theta)$ which include examples like $q(\theta)\,=\, (\theta+b)^{-k}$ for $k > 0$, $b \in \mathbb{R}$, and $\theta > -b$. \\

\noindent   Further applications include unbiased estimation: {\bf (i)} for truncated distributions generated from an exponential family density, such as truncated normal distributions, and {\bf (ii)} ratios of two exponential family natural parameters, including the ratio of normal means extending also to the general bivariate normal model.

\subsection{Motivation}
\label{motivation}
For normal models with known variance,  the exponential family natural parameter $\theta$ (see (\ref{modeldensity})) is a multiple of the mean, hence, the estimation of inverse moments, in particular with the criterion of unbiasedness, relates to estimating inverse powers of $\theta$.  The problem of making inferences about $\frac{1}{\theta}$ is thus of interest for normal models (see for instance \cite{withers2013inverse}, and the references therein) as well as for exponential families as in (\ref{modeldensity}).    Scientific and economic applications include:

\begin{itemize}
    \item Experimental Nuclear Physics: The momentum of a charged particle is determined by measuring its track curvature, a technique widely used in high-energy physics experiments (\cite{lamanna1981curvative}, \cite{treadwell1982momentum}). Estimating inverse moments of the curvature parameter allows for improved precision in momentum calculations.
    \item Optimal Control Theory: In dynamic decision-making problems, particularly in the one-dimensional special case of the single-period control problem, inverse moment estimation is critical for deriving optimal strategies (\cite{zaman1981reciprocal}, \cite{zellner1971lognormal}).
    \item Economic Modeling: Inverse moments frequently arise in macroeconomic and microeconomic models. For instance, the investment multiplier in Keynesian economics is a function of the marginal propensity to consume, while the long-run supply elasticity in agricultural economics can be estimated from Nerlove’s supply response model (\cite{braulke1982note}).
    \item Astronomical Data Analysis: In observational astronomy, data are often truncated due to detection limits—only sources above a brightness threshold are observed. When measurements follow an exponential distribution, as in photon counts or luminosity, inverse moments of the natural parameter (e.g., $1/\lambda$) are key for estimating population means and correcting selection bias. Unbiased estimators of these moments are thus essential for accurate inference under truncation (\cite{efron1999nonparametric}).

\end{itemize}
These applications highlight the widespread relevance of unbiased estimation of inverse moments in diverse scientific fields, motivating the need for a general and systematic approach to their construction whenever possible.

\section{Main results and applications}

\noindent   Our findings are applicable to statistics $X \in \mathbb{R}$ with exponential family densities 
\begin{equation}
\label{modeldensity}
X \sim \, \,  h(x) \, \exp\{x\theta - \mathcal{A} (\theta)\},
\end{equation}
with respect to Lebesgue measure on $(-\infty,a)$ where $a$ can be either finite or ``$+\infty$'', as well as to extensions through changes of (i) location, (ii) scale, (iii) sign, and (iv) truncation of $X$ to $(-\infty,b)$ with $b<a$.
We assume that the natural parameter space $\Theta=\{\theta \in  \mathbb{R}: \int_{-\infty}^a h(x) \,exp(x\theta)  \, dx < \infty  \}$ contains an interval of the form $(0, \zeta)$. 

\noindent  Unless otherwise specified, we consider reduced exponential family densities in (\ref{modeldensity}), but we point out the well-known property (e.g., \cite{efron2022exponential})
that such a model arises with $X_1, \ldots, X_n$ independently distributed and drawn from exponential family density  $g(u) \, \exp\{\theta t(u) - \mathcal{C}(\theta)  \}$, since the statistic $T=\sum_{i=1}^n t(x_i)$ is a complete and sufficient statistic for the family of densities with density representable as in (\ref{modeldensity}).    Moreover, completeness implies that there exists at most one unbiased estimator of a parametric function $\gamma(\theta)$ which is a function of $T$.   Before continuing, we point out that for Bayesian inference, conjugate prior distributions in one-parameter exponential families are readily available (e.g., \cite{diaconis1979conjugate}), maintain the exponential family structure, and facilitate the estimation of inverse moments.
\noindent 
We first present a general non-existence unbiased estimation result for inverse powers of
$\theta$.  
\begin{theorem}
\label{theoremnon-existence2}
Let $X$ have density in (\ref{modeldensity}) and let the interior of $\Theta$ include the point $\theta=0$.  Then there does not exist an unbiased estimator for $\frac{1}{\theta^k},\quad k >0$. 
\end{theorem}
\noindent  {\bf Proof.}
\noindent Suppose there exists an unbiased estimator, say $\delta(X)$, of $\frac{1}{\theta^k}$.  Then, we would have $$\mathbb{E} \{ \delta (X) \} \,  =  \, \int_{-\infty}^{\infty} \delta (x) \, h(x) \, \exp\{(x \theta - \mathcal{A} (\theta)\} \,dx \, =  \, \frac{1}{\theta^k}\,.$$
Since the left-hand side is analytic in $\theta$ (e.g., \cite{barndorff1978information}, \cite{brown1986fundamentals}), while the right hand side does not admit derivatives of all orders at $\theta = 0$, the equation cannot be satisfied in $\delta(\cdot)$, and therefore there does not exist an unbiased estimator for $\frac{1}{\theta^k}$.  \qed  \\

\noindent In contrast, the next result provides an unbiased estimator of $\frac{1}{\theta}$ under the condition that $\theta>0$.
\begin{proposition}
\label{proposition-main}
For $X$ with density as in (\ref{modeldensity}), there exists for $\theta>0$ a unique unbiased estimator of $\frac{1}{\theta}$ and it is given by
\begin{equation}
\label{unbiasedestimator}
\delta_0(X) \, = \, \frac{1}{h(X)} \, \int_X^{a}  h(s) \, ds\,.
\end{equation}
\end{proposition}
\noindent {\bf Proof.}
We have for all $\theta \in \Theta$
\begin{eqnarray*}
    \mathbb{E}_{\theta} \big( \delta_0 (X)\big) 
     \, & = & \,  \int_{-\infty}^{a} \int_{x}^{\infty} h(s) \, \exp\{(x\theta - \mathcal{A} (\theta)\} \, ds  \,dx \,  \\
       \, & = & \,  \int_{-\infty}^a h(s) \int_{-\infty}^{s} \, exp\{x\theta - \mathcal{A} (\theta)\} \,dx \,ds  \quad  (\hbox{Fubini's theorem})\\ \, & = & \,  \int_{-\infty}^{\infty}  h(s) \,\frac{1}{\theta}  \, exp\{s\theta - \mathcal{A} (\theta)\}  \,ds \quad  (\hbox{since } \theta \hbox{ is positive}) \\ 
       \, & = & \,  \frac{1}{\theta} \,. \quad \quad \quad \quad \quad \quad \quad \quad \quad \quad \quad \quad \quad  \quad \quad \quad \quad \quad \quad \quad \quad \quad \quad  \text{\qed}
\end{eqnarray*} 

\begin{remark}
Denoting the densities in (\ref{modeldensity}) as $f_{\theta}$ and their survival functions as $\bar{F}_{\theta}$, notice that $f_0(x)\,=\, h(x) \, e^{-\mathcal{A} (\theta)}$ so that the unbiased estimator in (\ref{unbiasedestimator}) is expressible as
\begin{equation}
\nonumber
\delta_0(X) \, = \, \frac{\bar{F}_0(X)}{f_0(X)}\,.
\end{equation}
   
\end{remark}
Before pursuing with various extensions, both with respect to model and parametric function $\gamma(\theta)$ variations, let us simply record Voinov's result (\cite{voinov1985unbiased}) applicable to a normal model with known variance.   Denote $\phi$ and $\bar{\Phi}$ as the probability density and survival functions of a $N(0,1)$ distribution. 

\begin{example}
Let $Y \sim N(\mu, \sigma^2)$ with $\mu>0$.   Then $X = \frac{Y}{\sigma} \sim N(\theta, 1)$ has model density (\ref{modeldensity}) with $\theta=\frac{\mu}{\sigma}$ and $h(x)=\phi(x)$.  Applying Proposition \ref{proposition-main}, we obtain that $\delta_0(X)\, = \, \frac{\bar{\Phi}(X)}{\phi(X)}$ is unbiased for $\frac{1}{\theta}$ or, equivalently, that 
$ \, \frac{1}{\sigma} \, \frac{\bar{\Phi}(\frac{Y}{\sigma})}{\phi(\frac{Y}{\sigma})}$ is unbiased for $\frac{1}{\mu}$. 
\end{example}
\begin{corollary}
\label{corollary-extensions}
\begin{enumerate}
\item[ (a)]   (change of location)    For density (\ref{modeldensity}) with $\theta > \theta_0$, the estimator   

\begin{equation}
\label{unbiasedestimatorlocation}
\delta_1(X) \, = \, \frac{1}{h(X)  \, e^{X \theta_0}} \, \int_X^{a}  h(s) \, \, e^{s\theta_0}  ds\, = \, \frac{\bar{F}_{\theta_0}(X)}{f_{\theta_0}(X)}
\end{equation}
is an unbiased estimator of $\frac{1}{\theta- \theta_0}$.

\item[ (b)] (change of sign)    For $X$ having density $h(x)\, \exp{\{\theta x - \mathcal{A}(\theta)\}}$ with respect to Lebesgue measure on $(a',\infty)$ and $\theta<0$, the estimator $\delta_2(X) \, = \, - \frac{F_0(X)}{f_0(X)}$ is unbiased for $\frac{1}{\theta}$.
Equivalently, for model density 
\begin{equation}
\label{modeldensity2}
X \sim h(x)\, \exp{\{-\theta' x - \mathcal{A}(-\theta')\}}
\end{equation}
 with respect to Lebesgue measure on $(a',\infty)$ and $\theta'>0$, the estimator $-\delta_2(X)\, = \,  \frac{F_0(X)}{f_0(X)}$ is unbiased for 
$\frac{1}{\theta'}$.

\item[ (c)]   (truncation)  For an observable $T$ with density 
\begin{equation}
\label{densitytruncated}  T \sim h(t) \, \mathbb{I}_{(-\infty,b)}(t) \, \exp\{\theta \, t \, - \, \mathcal{A}(\theta) \, - \log \mathbb{P}_{\theta}(X \leq b)   \},
\end{equation}
which is the truncated to $(-\infty,b)$ version of $X$ in (\ref{modeldensity}), the estimator
\begin{equation}
\delta_3(T) \, = \, \frac{1}{h(T)} \, \int_T^{b}  h(s) \, ds\,
\end{equation}
is an unbiased estimator of $\frac{1}{\theta}$.
\end{enumerate}

\end{corollary}
\noindent {\bf Proof.}
\begin{enumerate}
\item[ (a)]   The result follows by expressing the density for $X$ as $$g(x) \, \exp\{x \theta'\, - \, \mathcal{A}(\theta) \big\}\,,$$
with $g(x)\,=\, h(x) e^{x \theta_0}$ and $\theta'=\theta-\theta_0$, and then applying Proposition \ref{proposition-main} to $g$ and $\theta'>0$.

\item[ (b)]   Let $\theta'=-\theta>0$ and $a=-a'$.  Since $Y \overset{d}{=} - X$ has density $h(-y) \exp\{\theta' y - A(-\theta')  \}$ with respect to Lebesgue measure on $(-\infty, a)$, we can apply Proposition \ref{proposition-main} to infer that
\begin{equation}
\nonumber
- \delta_0(Y) \, = \, - \frac{\int_Y^{a} \, h(-s) \, ds }{h(-Y)} \, = \, - \frac{\int_{a'}^{-Y} \, h(s) \, ds }{h(-Y)} \,= \, \delta_2(X)
\end{equation}
is unbiased for $- \, \frac{1}{\theta'}\,=\, \frac{1}{\theta}$.
\item[ (c)] The result follows directly from Proposition \ref{proposition-main} with $h(t)$ replaced by $h(t) \, \mathbb{I}_{(-\infty,b)}(t)$.    \qed

\end{enumerate}

\noindent  Combinations of Proposition \ref{proposition-main} and different parts of Corollary \ref{corollary-extensions} produce a wealth of applications, and here is a selection of examples.

\begin{example}
(a) For $X \sim N(\theta, 1)$ with the lower bound parametric restriction $\theta \geq \theta_0$, part (a) of Corollary \ref{corollary-extensions} applies and yields that
$$\delta_1(X) \, = \, \frac{\bar{\Phi}(X-\theta_0)}{\phi(X-\theta_0)}$$
is unbiased for $\frac{1}{\theta-\theta_0}$.

\noindent (b)   For $X \sim N(\theta, 1)$ with the upper bound parametric restriction $\theta \leq 0$, part (b) of Corollary \ref{corollary-extensions} applies and yields that
$$ \delta_2(X) \, = \, - \frac{\Phi_0(X)}{\phi_0(X)}$$ is unbiased for $\frac{1}{\theta}$.  
\end{example}
\begin{example} (truncated normal distribution)
\label{exampletruncatednormal}
Proposition \ref{proposition-main} or part (c) of Corollary \ref{corollary-extensions} apply directly to truncated normal distributions such that $T \overset{d}{=} Z|Z \leq b$ with 
$Z \sim N(\mu, \sigma^2)$ and $\mu>0$.  Indeed $T$ has density given in (\ref{densitytruncated}) with $h(t)\,=\,\frac{1}{\sigma} \, \phi(\frac{t}{\sigma}) \,$ and $\theta\,=\, \frac{\mu}{\sigma^2}$.  Since $\theta>0$, we have that
$$\delta_3(T) \, = \, \frac{\sigma}{\phi(\frac{X}{\sigma})} \,  \int_T^b   \frac{1}{\sigma} \, \phi(\frac{s}{\sigma}) \, ds \, = \, \sigma \, \frac{\big\{ \bar{\Phi}(\frac{T}{\sigma}) \, - \, \bar{\Phi}(\frac{b}{\sigma})\big\}}{\phi(\frac{T}{\sigma})} $$
is an unbiased estimator of $\frac{1}{\theta}\,=\, \frac{\sigma^2}{\mu}$, 
and that $\frac{1}{\sigma^2} \, \delta_3(T)$ is therefore unbiased for $\frac{1}{\mu}$.
\end{example}
\begin{remark}
\label{remarknormalonly}
In Example \ref{exampletruncatednormal} where $\mathbb{E}(T)\, = \, \mu + \sigma \frac{\phi(\frac{b-\mu}{\sigma})}{\Phi(\frac{b-\mu}{\sigma})}$, we point out that the results of the paper do not address the estimation of  $\frac{1}{E(T)} \neq  \frac{1}{\mu}$, but rather that of  $\frac{1}{\mu}$.   More generally for model densities in (\ref{modeldensity}), we have $\mathbb{E}_{\theta}(X)\,=\, A'(\theta)$, so that $\mathbb{E}_{\theta}(X)\,= c_1 \theta$ iff $A(\theta)\,=\, \frac{c_1 \, \theta^2}{2} \, + \, c_2$, for some constants $c_1$ and $c_2$, which in turn occurs only for $X \sim N(\theta, \sigma^2)$
\end{remark}

\begin{example} (Gamma model)   Let $X \sim \mathcal{G}(\alpha, \theta')$ with $\alpha>0$ known $\theta'>0$, and density $\frac{\theta'^{\alpha}}{\Gamma(\alpha)} \, x^{\alpha-1}  \, e^{-\theta' x} \, \mathbb{I}_{(0,\infty)}(x)$.  
Obviously $\mathbb{E}(\frac{X}{\alpha}) \, = \, \frac{1}{\theta'}$, so our findings above are not required, but it is nevertheless interesting to illustrate that this unbiased estimator can be generated from part $(b)$ of Corollary \ref{corollary-extensions} with $h(x)\,=\, x^{\alpha-1}, a'=0$, and 
\begin{equation}
\nonumber   - \delta_2(X) \, = \,  \frac{1}{X^{\alpha-1}} \, \int_0^X   s^{\alpha-1} \, ds \, = \, \frac{X}{\alpha}.  
\end{equation}
Although the initial motivation for our paper, as well as for Voinov's work, involved a positivity restriction for a normal mean, observe that there is no parametric restriction for this Gamma model example, and that it is not required in Proposition \ref{proposition-main}. 

\noindent Finally, observe that the illustration extends nicely to a truncated Gamma version of $X$ on $(b,\infty)$ with $b>0$, i.e., for $T \overset{d}{=} X|X>b$.  Indeed, as in part (c) of Corollary \ref{corollary-extensions}, we can proceed as above to infer that
\begin{equation}
\nonumber   - \delta_3(T) \, = \,  \frac{1}{T^{\alpha-1}} \, \int_b^T   s^{\alpha-1} \, ds \, = \, \frac{T}{\alpha} \, - \, \frac{b^{\alpha}}{ \alpha \, T^{\alpha-1}}
\end{equation}
is an unbiased estimator of $\frac{1}{\theta'}$.
\end{example}

\begin{example}
\label{exampleinversegaussian}  (Inverse Gaussian model)
Inverse Gaussian models ($X \sim \hbox{IG}(\mu, \lambda)$) have densities which can be expressed as 
\begin{equation}
\nonumber
f_{\mu, \lambda}(x)\,=\,  \lambda^{1/2} \, (2\pi x^3)^{-1/2}  \,  \exp\big\{-\frac{\lambda x}{2 \mu^2} \,- \, \frac{\lambda}{2x}\,+ \, \frac{\lambda}{\mu}   \big\} \, \mathbb{I}_{(0,\infty)}(x)\,,
\end{equation}
with $\mu>0$ a mean parameter and $\lambda>0$ a shape parameter.   Now, consider $\lambda$ to be known in which case $X$ has density expressible as in (\ref{modeldensity2}) with $a'=0$, $\theta'=\frac{\lambda}{2 \mu^2}>0$, and $h(x)\,=\, x^{-3/2}\, e^{-\frac{\lambda}{2x}}  \, \mathbb{I}_{(0,\infty)}(x)\, $.   Part (b) of Corollary \ref{corollary-extensions} then applies and an evaluation which exploits the presence of a $\chi^2_1$ density (or an incomplete Gamma function) tells us that
\begin{eqnarray*}
- \delta_2(X) \, & = & \frac{\int_0^X  \, s^{-3/2} \, e^{-\lambda/2s} \, ds }{X^{-3/2} \, e^{-\lambda/2X}} \\
\, & = & (\lambda \, X^3)^{1/2}  \, e^{\lambda/2X} \, \int_{\lambda/X}^{\infty}  \, t^{-1/2} \, e^{-t /2} \, dt \\
\, & = &  \,  \big(\frac{8 \pi X^3}{\lambda}\big)^{1/2} \, e^{\lambda/2X} \, \bar{\Phi}\big(\sqrt{\frac{\lambda}{X}}  \big)
\end{eqnarray*}
is an unbiased estimator of $\frac{1}{\theta'}$.  Observe that, in terms of the original parameters, the estimator  $\frac{\lambda \delta_2(X)}{2}$ is unbiased for $\mu^2$.    To conclude, we point out that unbiased estimators of various functions of $(\mu,\lambda)$ have been proposed in the literature, namely by \cite{iwase1983umvue}, but such results to the best of our knowledge do not cover the example here (also see \cite{seshadri1993inverse}, chapter 6).
\end{example}

We conclude this section with an application to estimating unbiasedly the ratio of two normal means based on independent samples, as well as based on a bivariate normal distribution.
\begin{example}  (ratio of normal means)
\begin{enumerate}
\item[ (A)]   Consider $Z_i,  i=1,  \ldots, n_1+n_2,$ independently distributed as 
$N(\mu_1, \tau_1^2)$ for $Z_1, \ldots, Z_{n_1}$, and $N(\mu_2, \tau_2^2)$ for $Z_{n_1+1}, \ldots, Z_{n_1+n_2}$, with $\tau_1^2$ and $\tau_2^2$ known.  With $\bar{Z}_1=\frac{1}{n_1}\sum_{i=1}^{n_1} Z_i$ and $X \,=\, \frac{\sqrt{n_2}}{\tau_2}\bar{Z}_2\, = \, \frac{1}{\sqrt{n_2} \, \tau_2}\sum_{i=n_1+1}^{n_1+n_2} Z_i$ jointly sufficient for $(\mu_1, \mu_2)$, with $X$ having density as in (\ref{modeldensity}) with $a=\infty$ and $\theta\,=\, \frac{\sqrt{n_2} \, \mu_2}{\tau_2}$, it follows from Proposition \ref{proposition-main}, as in Example $(1)$, where $\delta_0(X)\,=\, \frac{\bar{\Phi}(X)}{\phi(X)}$ that $\mathbb{E}\big(\delta_{0}(X)\big) \, = \, \frac{\tau_2}{\sqrt{n_2}\, \mu_2}$ 
 for all $\mu_2>0$, and that the estimator $\frac{\sqrt{n_2}}{\tau_2}\, \bar{Z}_1 \,  \delta_{0}(\frac{\sqrt{n_2}}{\tau_2}\bar{Z}_2)$ is an unbiased estimator of $\frac{\mu_1}{\mu_2}$ under the restriction $\mu_2>0$ and given the independence.

\item[ (B)]  The result in part (A) also admits an extension to the bivariate normal case with non-independent components.  Without loss of generality, consider $(Y_1, Y_2)^{\top} \sim N_2(\mu, \Sigma)$ with  $\mu=(\mu_1, \mu_2)^{\top}$, \[\Sigma=\left[\begin{array}{cc}
\sigma_1^2 & \rho \sigma_1 \sigma_2   \\
\rho \sigma_1 \sigma_2 & \sigma_2^2    
\end{array}  \right], \]
and known $\rho, \sigma_1, \sigma_2$. \footnote{ Normally, we would expect to require that $|\rho|<1$ but the result still holds here for $|\rho|=1$ with the particularity that $W$ is degenerate.}   Indeed since $W=Y_1 - \rho \frac{\sigma_1}{\sigma_2} Y_2$ and $Y_2$ are independently distributed as $N(\mu_1 - \rho \frac{\sigma_1}{\sigma_2} \mu_2, 1- \rho^2)$ and $N(\mu_2, \sigma_2^2)$, respectively, we obtain from Proposition \ref{proposition-main} and the independence property that:
\begin{equation}
\nonumber   \mathbb{E}\Big(W \, \frac{1}{\sigma_2}\delta_{0}\left(\frac{Y_2}{\sigma_2}\right)  \Big)\, = \, \big( \mu_1 - \rho \frac{\sigma_1}{\sigma_2} \mu_2 \big) \, \frac{1}{\mu_2} \, = \,  \frac{\mu_1}{\mu_2} - \rho \frac{\sigma_1}{\sigma_2} \, \hbox{ for } \mu_2>0, 
\end{equation}
which implies that $W  \frac{1}{\sigma_2}\delta_{0}\left(\frac{Y_2}{\sigma_2} \right)+ \rho \frac{\sigma_1}{\sigma_2}$ is unbiased 
for $\frac{\mu_1}{\mu_2}$, whenever $\mu_1 \in \mathbb{R}$ and $\mu_2>0$. 
\end{enumerate}
\end{example}

\subsection{Extensions to a class of parametric functions}

\noindent The next result generalizes the unbiased estimation finding of Proposition \ref{proposition-main} to parametric functions that are expressible as
\begin{equation}
\label{q}
q(\theta)  \, = \, \int_0^{\infty} f(y) \, \exp(- y \theta) \, dy \,,
\end{equation}
for a given integrable function $f$. Unbiasedness will be established for values of $\theta$ such that $q(\theta)$, i.e., for $\theta \in \mathcal{C}=\{\theta \in  \mathbb{R}: |q(\theta)|<\infty    \}$.   The functions $q(\theta)$ are precisely the class of functions that are completely monotone, which are also characterized by functions $q$ with derivatives of alternating signs,
i.e.,  $ (-1)^m \, \frac{d^m}{d \theta^m}  \, q(\theta)  \geq 0$ for all $\theta \in \mathcal{C}$, and for $m=1,2, \ldots$, (e.g., \cite{feller1971introduction}, \cite{widder1941laplace}).      Immediate applications of interest generated by (\ref{q}) 
include: (i) the basic case $q(\theta) \, = \, \frac{1}{\theta}$ for $f(y)=1$ and $\mathcal{C}=(0,\infty)$, (ii) the power extensions $q(\theta) \, = \, \frac{1}{\theta^k}, k >0 ,$ for $f(y)\,=\, \frac{1}{\Gamma(k)} y^{k-1}$ and $\mathcal{C}=(0,\infty)$, (iii) further extensions $q(\theta) \, = \, \frac{1}{(b+\theta)^k}$, $k >0$, for $f(y)\,=\, \frac{1}{\Gamma(k)} y^{k-1} e^{-by}$ and $\mathcal{C}=(-b,\infty)$, and many more as there few restrictions on $f$, including for instance $q(\theta)\,=\, \frac{e^{-d_1\theta}-e^{-d_2\theta}}{\theta}$ generated by the indicator function $f(y)\,=\, \mathbb{I}_{(d_1, d_2)}(y)$ with $0 \leq d_1 < d_2 \leq \infty$.
\begin{proposition}
\label{proposition-main-extension}
For $X$ with density as in (\ref{modeldensity}), there exists for $\theta \in \mathcal{C}$ a unique unbiased estimator of $q(\theta)$ of the form (\ref{q}),  and it is given by
\begin{equation}
\label{delta0f}
    \delta_0(X) \, = \,  \frac{1}{h(X)} \, \int_{X}^{a} \, h(s) \, f(s-X)\,  ds 
\end{equation}
{\bf Proof.  }  We have for $\theta \in \mathcal{C}$
\begin{eqnarray*}
 \mathbb{E}_{\theta} \{\delta_0(X)\}   \,&\, = &  \int_{-\infty}^{a} \int_{x}^{a} \, h(s) \,f(s-x)  \,ds\, \exp\{x \, \theta - \mathcal{A} (\theta)\} \,dx \,  \\ 
\, & \, = & \int_{-\infty}^{a} \, h(s) \,  \int_{-\infty}^{s} \,f(s-x) \, \exp\{x \, \theta - \mathcal{A} (\theta)\} \,dx\, ds \, (\hbox{ Fubini theorem})\\
\, & \, = & \int_{-\infty}^{a} \, h(s) \,  \int_{0}^{\infty} \, \exp\{s\, \theta - u \theta - \mathcal{A} (\theta)\} \, f(u) \, du\, ds \, \\
& & \quad \quad \quad \quad \quad \quad \quad \quad (\hbox{change of variable } u=x-s)\\ 
&= &   \, q(\theta) \, \int_{-\infty}^{a}  h(s) \exp\{s \theta - \mathcal{A} (\theta)) \,  ds   \\
\, & = & \, q(\theta)\,. \quad \quad \quad \quad \quad \quad \quad \quad \quad \quad \quad \quad  \quad \quad \quad \quad \quad \quad \quad \quad \quad \quad \quad \quad \quad \text{\qed}
\end{eqnarray*}
\end{proposition}

A vast number of examples follow from Proposition \ref{proposition-main-extension} for different choices of $h$ in (\ref{modeldensity}) and different $f$'s in (\ref{q}) or (\ref{delta0f}), and then with changes of sign or truncation as in Corollary \ref{corollary-extensions}.     Here are such examples.

\begin{example}  
For $X \sim N(\theta, 1)$ and $q(\theta)\,=\, (b+\theta)^{-k}$, $b \in \mathbb{R}, k >0$, Proposition \ref{proposition-main-extension} applies for $a=\infty$, $h(t)=\phi(t)$, $f(y)\, = \, \frac{1}{\Gamma(k)} \, y^{k-1} e^{-by}$ and $\mathcal{C}=(-b,\infty)$, yielding with a change of variables the unbiased estimator 
\begin{eqnarray*}
  \delta_0(X) \,  & = & \,  \frac{1}{\Gamma(k) \, \phi(X)}   \int_X^{\infty} \phi(s) \, (s-X)^{k-1} \, e^{- b(s-X)} \, ds \\
\, & = & \,    \frac{1}{\Gamma(k) \, \phi(X+b)} \,   \int_{X+b}^{\infty} \phi(u) \, (u-b-X)^{k-1} \, du
\end{eqnarray*}
for $q(\theta), \theta>-b$.
\end{example}
From this, Example $2$ (a) is recovered with $k=1$ and $\theta_0=-b$.   For $k=2$, one obtains with the evaluation $\int\limits_{X + b}^{\infty}\;u\phi(u)\;du = \phi(X + b)$, then
\begin{equation}
\nonumber \delta_0(X) \, = \,  \big\{1 \, - \, (X+b) \, \frac{\bar{\Phi}(X+b)}{\phi(X+b)}                  \big\}
\end{equation}
is an unbiased estimator of $(\theta+b)^{-2}$ for $\theta>-b$, $b \in  \mathcal{R}$.

\section*{Concluding remarks}
Our findings shed light on the existence of an unbiased estimator $\delta_0(X)$ of $\frac{1}{\theta}$ for $X \sim N(\theta, \delta^2)$ with $\theta>0$, by inscribing the result in a more general finding, applicable to: (i) the class of one-parameter continuous exponential families with $\theta$ being the natural parameter, and (ii) to parametric functions $q(\theta)$ which are completely monotone.  

 A related issue, necessarily specific to the given context (i) and (ii), concerns the point estimation performance of the estimator $\delta_0(X)$ of $q(\theta)$ in comparison to other estimators, other than unbiasedness which follows from our results.   Such comparisons require caution and the possible use of a bounded loss function or a sufficiently slow varying one.  Indeed, with the routine choice of mean squared error for the normal case with $\sigma^2=1$ for instance, we have using Remark 1 for all $ \theta>0$:
 \begin{equation}
 \nonumber
 \hbox{MSE}(\theta, \delta_0) \, = \, \mathbb{E} \big(\frac{\bar{\Phi}^2(X)}{\phi^2(X)}   \big) \, = \, \int_{\mathbb{R}} g_{\theta}(x)\,dx \, - \frac{1}{\theta^2},
 \end{equation}
with $g_{\theta}(x)\,=\, \sqrt{2 \pi} \, e^{-\theta^2/2} \, \bar{\Phi}^2(x) \, \exp\{\theta x + x^2/2\}$.   Since, for all $\theta>0$, $\lim_{x \to \infty} g_{\theta}(x) \,=\, + \infty$, the mean squared error of $\delta_0$ does not exist and cannot be used directly for evaluation purposes.


\end{document}